\documentclass[11pt, twoside]{amsart}

\usepackage[utf8]{inputenc}
\usepackage[T1]{fontenc}
\usepackage{amssymb,amsmath,amstext, graphicx, hyperref}

\newtheorem{theor}{Theorem}[section]
\newtheorem{lem}[theor]{Lemma}
\newtheorem{defin}[theor]{Definition}

\newtheorem{prop}[theor]{Proposition}

\newtheorem{cor}[theor]{Corollary}
\newtheorem{rem}[theor]{Remark}

\newtheorem{terminologyandnotation}[theor]{Terminology and notation}

\numberwithin{equation}{section}

\newcommand{\lan}{\langle}
\newcommand{\ran}{\rangle}

\newcommand{\es}{\emptyset}

\newcommand{\uhr}{\upharpoonright}

\newcommand{\nts}{\negthickspace}
\newcommand{\uhrc}{\nts \upharpoonright \nts}

\newcommand{\mcA}{\mathcal{A}}
\newcommand{\mcB}{\mathcal{B}}

\newcommand{\mcG}{\mathcal{G}}

\newcommand{\mcM}{\mathcal{M}}
\newcommand{\mcN}{\mathcal{N}}

\newcommand{\mbS}{\mathbf{S}}

\newcommand{\mbbN}{\mathbb{N}}

\newcommand{\mbbZ}{\mathbb{Z}}
\newcommand{\mbbQ}{\mathbb{Q}}

\newcommand{\Aut}{\mathrm{Aut}}

\newcommand{\Spt}{\mathrm{Spt}}
\newcommand{\spt}{\mathrm{spt}}

\title[Typical automorphism groups]
{Typical automorphism groups of \\ finite nonrigid structures}

\author{Vera Koponen}

\address{Vera Koponen, Department of Mathematics, Uppsala University, Box 480,
75106 Uppsala, Sweden.}

\email{vera@math.uu.se}

\begin{document}

\begin{abstract}
We work with a finite relational vocabulary with at least one relation symbol with arity at least 2.
Fix any integer $m > 1$. For almost all finite structures (labelled or unlabelled)  such that at least $m$ elements are moved by some
automorphisms, the automorphism group is $(\mbbZ_2)^i$ for some $i \leq (m+1)/2$; and if some relation symbol has
arity at least 3, then the automorphism group is almost always $\mbbZ_2$.

\end{abstract}

\maketitle

\section{Introduction}

\noindent
This article complements the work in \cite{AK} with quite explicit information about the automorphism group of
``almost all'' finite structures such that at least $m$ elements are moved by some automorphisms, for any fixed integer $m$.
It turns out that the automorphism group is almost always a power of $\mbbZ_2$, where the maximal power is bounded by $(m+1)/2$.
As part of proving this we prove that almost all finite structures such that at least $m$ elements are moved by some automorphisms have the
property that exactly $m'$ elements are moved by some automorphism, where $m' = m$ if $m$ is even and $m' = m+1$ otherwise.
Perhaps surprisingly, we get different results depending on the maximal arity of the relation symbols (of the finite relational language).
If the maximal arity is at least 3, then the typical automorphism group is always $\mbbZ_2$, no matter what $m$ is.
If the maximal arity is 2, then for each $i = 1, \ldots, m'/2$, $(\mbbZ_2)^i$ appears as an automorphism group with positive probability
(given by the uniform probability measure on the set of $n$-element structures).
The situation is slightly different if we restrict attention to finite structures such that {\em exactly} 
$m$ elements are moved by some automorphisms.
Then $\mbbZ_3$ or the symmetric group on three elements appear as a subgroup of the typical automorphism group if $m$ is odd.
These results hold for both labelled and unlabelled structures (See Remark~\ref{remark about unlabelled structures}).

We now introduce some notation and terminology which will be used throughout the article and then
state the two main results.
We work with a finite relational vocabulary (also called signature) $\{R_1, \ldots, R_\rho\}$,
where each relation symbol $R_i$ has arity $r_i$.
The number $r = \max\{r_1, \ldots, r_\rho\}$ is called the {\em maximal arity} and the we assume that it is at least 2.
The set of all structures for this vocabulary with universe $[n] = \{1, \ldots, n\}$ is denoted $\mbS_n$ and we let
$\mbS = \bigcup_{n=1}^{\infty}\mbS_n$.
For any set $A$, $|A|$ is its cardinality and $Sym(A)$ the group of all permutations of $A$.
Suppose that $f_1, \ldots f_k \in Sym(A)$.
Then $\lan f_1, \ldots, f_k \ran$ denotes the subgroup of $Sym(A)$ generated by $f_1, \ldots f_k$ and we define
\[ 
\Spt(f_1, \ldots, f_k) \ = \ \{a \in A : g(a) \neq a \text{ for some } g \in \langle f_1, \ldots, f_k \rangle \}
\]
and let $\spt(f_1, \ldots, f_k) = |\Spt(f_1, \ldots, f_k)|$.
We call $\Spt(f_1, \ldots, f_k)$ the {\em support} of $f_1, \ldots, f_k$.
For any finite structure $\mcM$ we let $\Aut(\mcM)$ denote its group of automorphisms,
\begin{align*}
\spt(\mcM) \ &= \ \max\{\spt(f) : f \in \Aut(\mcM)\},\\
\Spt^*(\mcM) \ &= \ \{a \in M : a \in \Spt(f) \text{ for some } f \in \Aut(\mcM)\}, \ \text{ and}\\
\spt^*(\mcM) \ &= \ \big|\Spt^*(\mcM)\big|.
\end{align*}
We call $\Spt^*(\mcM)$ the {\em support} of $\mcM$. 
For every $m \in \mbbN$ define
\begin{align*}
\mbS_n(\spt \geq m) \ &= \ \{\mcM \in \mbS_n : \spt(\mcM) \geq m\} \ \ \text{ and} \\
\mbS_n(\spt^* \geq m) \ &= \ \{\mcM \in \mbS_n : \spt^*(\mcM) \geq m \}.
\end{align*}
Whenever $\mbS'_n \subseteq \mbS_n$ is defined for all $n \in \mbbN^+$ we let 
$\mbS' = \bigcup_{n=1}^{\infty} \mbS'_n$.
With the expression {\em almost all $\mcM \in \mbS'$ has the property $P$} we mean that
\[
\lim_{n\to\infty} \ \frac{\big| \{ \mcM \in \mbS'_n : \text{ $\mcM$ has $P$} \}\big|}{\big| \mbS'_n \big|}
\ = \ 1.
\]

\begin{theor}\label{characterisation of automorphism groups when r=2}
Suppose that the maximal arity is 2.
Let $m \geq 2$ be an integer and let $m' = m$ if $m$ is even and $m' = m+1$ otherwise.\\
(i) For almost all $\mcM \in \mbS(\spt \geq m)$, $\spt^*(\mcM) = m'$ and $\Aut(\mcM) \cong \big(\mbbZ_2\big)^i$
for some $i \in \{1, \ldots, m'/2\}$.\\
(ii) For every $i \in \{1, \ldots, m'/2 \}$ there is a rational number
$0 < a_i \leq 1$ (where $a_i < 1$ if $m > 2$) such that 
the proportion of $\mcM \in \mbS_n(\spt \geq m)$ such that $\Aut(\mcM) \cong \big(\mbbZ_2\big)^i$
converges to $a_i$ as $n \to \infty$.\\
(iii) Parts (i) and (ii) hold if `$\spt \geq m$' is replaced with `$\spt^* \geq m$'.
\end{theor}

\begin{theor}\label{characterisation of automorphism groups when r>2}
Suppose that $r$ is the maximal arity, $r \geq 3$ and let $m \geq 2$ be an integer.
Let $m' = m$ if $m$ is even and $m' = m + 1$ otherwise.
Then, for almost all $\mcM \in \mbS(\spt \geq m)$, $\spt^*(\mcM) = m'$ and $\Aut(\mcM) \cong \mbbZ_2$.
The same is true if `$\spt \geq m$' is replaced with `$\spt^* \geq m$'.
\end{theor}

\noindent
Intuitively, one may interpret the theorems as saying that finite structures tend to be as ``rigid'' as we allow them to be;
their automorphisms typically move as few elements as possible (given that impose the restriction that some minimum number
of elements are moved) and the automorphism group typically acts in the simplest possible way on the elements which are moved.
This is a generalisation of the well known result, proved via a sequence of articles \cite{ER, Fag77, FU, Har, Ober},
that almost all $\mcM \in \mbS$ are {\em rigid}, that is, $\Aut(\mcM)$ is trivial (i.e. contains only one element).

\begin{rem}\label{remark about unlabelled structures}{\rm
(i) Theorems~\ref{characterisation of automorphism groups when r=2}
and~\ref{characterisation of automorphism groups when r>2}
also hold if we consider {\em unlabelled} structures, that is, if we count structures only up to isomorphism.
This follows from the proof of Theorem~7.7 in \cite{AK}. \\
(ii) Theorems~\ref{characterisation of automorphism groups when r=2}
and~\ref{characterisation of automorphism groups when r>2} also hold if we require that all relations
are {\em irreflexive} or that all relations are {\em irreflexive and symmetric}, in the sense explained in
Remark~1.5 in \cite{AK}. Only minor modifications of the proofs (and some technical results) in \cite{AK} and this article are necessary.
}\end{rem}

\section{Preliminaries}\label{Preliminaries}

\begin{terminologyandnotation}{\rm
Recall the terminology and notation introduced before Theorem~\ref{characterisation of automorphism groups when r=2}.
So in particular we have fixed a finite relational vocabulary with maximal arity at least 2.
Structures (for this vocabulary) are denoted $\mcA, \mcB, \ldots, \mcM, \mcN$ and their universes
$A, B, \ldots, M, N$. For any set $A$, $|A|$ denotes its cardinality.
Since we mainly deal with structures $\mcM \in \mbS$, the universe will usually be $[n] = \{1, \ldots, n\}$ for some 
integer  $n > 0$.
For structures $\mcM$ and $\mcN$, $\mcM \cong \mcN$ means that they are isomorphic.
(See for example \cite{EF, Rot} for basic model theory.)
For groups $G$ and $H$, $G \cong H$ means that they are isomorphic as abstract groups.

Suppose that $f$ is a permutation of a set $\Omega$ and that $H$ is a group of permutations of $\Omega$.
Then $a \in \Omega$ is called a {\em fixed point of $f$} if $f(a) = a$. If $a$ is a fixed point of every $h \in H$,
then we say that $a$ is a {\em fixed point of $H$}.
For a structure $\mcA$ and $a \in A$, we call $a$ a {\em fixed point of $\mcA$} if $a$ is a fixed point of $\Aut(\mcA)$,
where we recall that $\Aut(\mcA)$ is the automorphism group of $\mcA$.
$Sym(\Omega)$ denotes the group of all permutations of $\Omega$, i.e. the symmetric group of $\Omega$,
and we let $Sym_n = Sym([n])$.

For a function $f A \to B$ and $X \subseteq A$, $f \uhrc X$ is the restriction of $f$ to $X$.
If $H$ is a permutation group on $\Omega$ and $X \subseteq \Omega$ is a union of orbits of $H$ on $\Omega$,
then $H \uhrc X = \{ h \uhrc X : h \in H\}$ and note that $H \uhrc X$ is a permutation group on $X$.
(For basic permutation group theory see \cite{DM} for example.)

It will be convenient to extend the notation used in the main results as follows:
\begin{align*}
\mbS_n(\spt^* = m) \ &= \ \{\mcM \in \mbS_n : \spt^*(\mcM) = m \},\\
\mbS_n(\spt^* \leq m) \ &= \ \{\mcM \in \mbS_n : \spt^*(\mcM) \leq m \},\\
\mbS_n(m \leq \spt^* \leq m') \ &= \ \{\mcM \in \mbS_n : m \leq \spt(\mcM) \leq m' \}.
\end{align*}
}\end{terminologyandnotation}

\noindent
We will use a some notions and results from \cite{AK} which we now state. 
The first gives an upper bound for $\spt^*(\mcM)$ for almost all $\mcM \in \mbS(\spt \geq m)$
and almost all $\mcM \in \mbS(\spt^* \geq m)$.

\begin{prop}\label{bound on support}\cite{AK}
For every $m \in \mbbN$ there is $m' \in \mbbN$ such that 
\begin{align*}
&\lim_{n\to\infty} \frac{\big| \mbS_n(\spt \geq m) \ \cap \ \mbS_n(\spt^* \leq m') \big|}
{\big| \mbS_n(\spt \geq m) \big|} \ = \\ 
&\lim_{n\to\infty} \frac{\big| \mbS_n(\spt^* \geq m) \ \cap \ \mbS_n(\spt^* \leq m') \big|}
{\big| \mbS_n(\spt^* \geq m) \big|}
\ = \ 1.
\end{align*}
\end{prop}

\noindent
Note that for every structure $\mcM$, $\Spt^*(\mcM)$ is the union of all nonsingleton 
orbits of $\Aut(\mcM)$ on $M$, 
so it makes sense to speak about $\Aut(\mcM) \uhrc \Spt^*(\mcM)$ and we always have
$\Aut(\mcM) \uhrc \Spt^*(\mcM) \cong \Aut(\mcM)$.

\begin{defin}\label{definition of S-n(A)}{\rm
Let $\mcA \in \mbS$ be such that $\Aut(\mcA)$ has no fixed point.
Suppose that $H$ is a subgroup of $\Aut(\mcA)$ such that $H$ has no fixed point.
For each integer $n > 0$, $\mbS_n(\mcA, H)$ is the set of $\mcM \in \mbS_n$ 
such that there is an embedding $f : \mcA \to \mcM$ such that $\Spt^*(\mcM)$ is the
image of $f$ and 
$H_f = \{f\sigma f^{-1} : \sigma \in H\}$ is a subgroup of $\Aut(\mcM) \uhrc \Spt^*(\mcM)$.
}\end{defin}

\begin{lem}\label{spt-star equals a union of S(A,H)}\cite{AK}
Let $m \geq 2$ be an integer. 
There are $\mcA_1, \ldots, \mcA_l \in \mbS_m$ without any fixed point and, for each 
$i = 1, \ldots, l$, subgroups $H_{i,1}, \ldots, H_{i,l_i} \subseteq \Aut(\mcA_i)$ 
without any fixed point such that
\[
\mbS(\spt^* = m) \ = \ \bigcup_{i=1}^l \bigcup_{j=1}^{l_i} \mbS(\mcA_i, H_{i,j}).
\]
\end{lem}

\noindent
Suppose that $H$ is a permutation group on $A$. 
Note that, for any integer $k >0$,  $H$ acts on $A^k$ by $h(a_1, \ldots, h_k) = (h(a_1), \ldots, h(a_k))$ for
every $(a_1, \ldots, a_k) \in A^k$ and every $h \in H$. Therefore we can speak of {\em orbits} of (the action of)
$H$ on $A^k$ for every $k > 0$.

\begin{defin}\label{definition of p(H)}{\rm
Suppose that $H$ is a permutation group on $A$.
Then $p(H) = |A|$, $q(H)$ is the number of orbits of $H$ on $A$ and $s(H)$ is the number of orbits of $H$ on  $A^2$.
 }\end{defin}

\begin{defin}\label{definition of beta}{\rm
Let
\[\beta(x,y,z) \ = \ k\binom{r}{2}x^2 \ - \ kr(r-1)xy \ - \ l(r-1)x \ + \ l(r-1)y \ + \ k\binom{r}{2}z.\]
}\end{defin}

\begin{prop}\label{comparisson of structures with different support-substructures}\cite{AK}
Suppose that $\mcA, \mcA' \in \mbS$ are such that neither $\Aut(\mcA)$ nor $\Aut(\mcA')$ has a
fixed point. Moreover, suppose that $H$ is a subgroup of $\Aut(\mcA)$ without fixed any point and 
that $H'$ is a subgroup of $\Aut(\mcA')$ without any fixed point. 
Let $p = p(H)$, $q = q(H)$, $s = s(H)$,
$p' = p(H')$, $q' = q(H')$ and $s' = s(H')$.\\
(i) The following limit exists in $\mbbQ \cup \{\infty\}$:
\[
\lim_{n\to\infty} \frac{\big|\mbS_n(\mcA', H')\big|}{\big|\mbS_n(\mcA, H)\big|}. 
\]
(ii) Suppose that $r = 2$.
\begin{itemize}
\item[(a)] If $p - q < p' - q'$ or if $p - q = p' - q'$ and $p > p'$, then
\[
\lim_{n\to\infty} \frac{\big|\mbS_n(\mcA', H')\big|}{\big|\mbS_n(\mcA, H)\big|} \ = \ 0.
\]
\item[(b)] If $p - q = p - q'$ and $p = p'$ then there is a rational number $a > 0$,
depending only on $\mcA$, $\mcA'$, $H$, $H'$ and the vocabulary, such that
\[
\lim_{n\to\infty} \frac{\big|\mbS_n(\mcA', H')\big|}{\big|\mbS_n(\mcA, H)\big|} \ = \ a.
\]
\end{itemize}
(iii) Suppose that $r > 2$ and let $\beta(x,y,z)$ be as in 
Definition~\ref{definition of beta}.
If any one of the two conditions
\begin{itemize}
\item[] $p - q < p' - q'$, or
\item[] $p - q = p' - q'$ and $\beta(p,q,s) > \beta(p',q',s')$
\end{itemize}
hold, then
\[
\lim_{n\to\infty} \frac{\big|\mbS_n(\mcA', H')\big|}{\big|\mbS_n(\mcA, H)\big|} \ = \ 0.
\]
\end{prop}

\begin{prop}\label{general quotient proposition}\cite{AK}
Let $\mcA_1, \ldots, \mcA_m, \mcA'_1, \ldots, \mcA'_{m'} \in \mbS$ 
be such that none of them has any fixed point.
Suppose that for every $i = 1, \ldots, m$ and $j = 1, \ldots, l_i$, 
$H_{i,j}$ is a subgroup of $\Aut(\mcA_i)$ without any fixed point and that
for every $i = 1, \ldots, m'$ and $j = 1, \ldots, l'_i$
$H'_{i,j}$ is a subgroup of $\Aut(\mcA'_i)$ without any fixed point.
Then the following limit exists in $\mbbQ \cup \{\infty\}$:
\begin{equation}\label{general limit of quotients of big unions}
\lim_{n\to\infty} \frac{\Big| \bigcup_{i=1}^{m'} \bigcup_{j=1}^{l'_i} \mbS_n(\mcA'_i, H'_{i,j}) \Big|}
{\Big| \bigcup_{i=1}^{m} \bigcup_{j=1}^{l_i} \mbS_n(\mcA_i, H_{i,j}) \Big|}.
\end{equation}
\end{prop}

\begin{defin}\label{definition of full permutation group}{\rm
Suppose that $\mcA \in \mbS$ has no fixed point and that $H$ is a 
subgroup of $\Aut(\mcA)$ without any fixed point.
For $\mcM \in \mbS_n(\mcA, H)$ we say that $H$ is the {\em full automorphism group of $\mcM$}
if for every isomorphism $f : \mcA \to \mcM \uhrc \Spt^*(\mcM)$ such
that $H_f = \{f \sigma f^{-1} : \sigma \in H\}$ is a subgroup of $\Aut(\mcM) \uhrc \Spt^*(\mcM)$
we have $H_f = \Aut(\mcM) \uhrc \Spt^*(\mcM)$.
}\end{defin}

\begin{lem}\label{automorhism group almost surely has same orbits as H}\cite{AK}
Suppose that $\mcA \in \mbS$ has no fixed point and that $H$ is a 
subgroup of $\Aut(\mcA)$ without any fixed point.\\
(i) For almost every $\mcM \in \mbS(\mcA, H)$,
$\Aut(\mcM) \uhrc \Spt^*(\mcM)$ has the same number of orbits as $H$.\\
(ii) Let  $G \leq H$. The proportion of $\mcM \in \mbS_n(\mcA, H)$ such that $G \cong \Aut(\mcM)$
converges to either 0 or 1 as $n \to  \infty$.
\end{lem}

\begin{lem}\label{typical number of nontrivial orbits of structures with bounded support}
Let $i$ be a positive integer.\\
(i) For almost all $\mcM \in \mbS(\spt^* = 2i)$, $\Aut(\mcM) \uhrc \Spt^*(\mcM)$
has exactly $i$ orbits on $\Spt^*(\mcM)$, so every such orbit has cardinality 2.\\
(ii) For almost all $\mcM \in \mbS(\spt^* = 2i+1)$, $\Aut(\mcM) \uhrc \Spt^*(\mcM)$
has exactly $i$ orbits, so $i-1$ orbits have cardinality 2 and the remaining orbit has cardinality 3.
\end{lem}

\noindent
{\bf Proof.}
We will use parts~(ii) and~(iii) of
Proposition~\ref{comparisson of structures with different support-substructures}.

(i) By Lemma~\ref{spt-star equals a union of S(A,H)},
there are $\mcA_1, \ldots, \mcA_m \in \mbS_{2i}$ without fixed points and for each
$i = 1, \ldots, m$ a number $l_i$ and subgroups $H_{i,1}, \ldots, H_{i,l_i}$ of $\Aut(\mcA_i)$
without fixed points
such that 
$\mbS_n(\spt^* = 2i) \ = \ \bigcup_{i=1}^m \bigcup_{j=1}^{l_i} \mbS_n(\mcA_i,H_{i,j})$ for each 
large enough $n$.
Moreover, 
by Lemma~\ref{automorhism group almost surely has same orbits as H},
for almost every $\mcM \in \mbS(\mcA_i, H_{i,j})$ the number of orbits of $H_{i,j}$ on $A_{i,j}$
is $q(H_{i,j})$.
Therefore it suffices to prove that there are $\mcA \in \mbS_{2i}$ without fixed point
and a subgroup $H \subseteq \Aut(\mcA)$ with
exactly $i$ orbits of cardinality 2 (then $H$ has no fixed points)
and that if $\mcA' \in \mbS_{2i}$ has no fixed point and
$H'$ is a subgroup of $\Aut(\mcA')$ without fixed points such that $H'$ does {\em not} have
exactly $i$ orbits of cardinality 2, then 
\begin{equation}\label{the proportion not having exactly i orbits of cardinality 2 approaches 0}
\lim_{n\to\infty} \frac{\big| \mbS_n(\mcA',H') \big|}{\big| \mbS_n(\mcA, H) \big|} \ = \ 0.
\end{equation}
First suppose that $\mcA \in \mbS_{2i}$ and that $H \subseteq \Aut(\mcA)$ has exactly $i$ orbits of
cardinality 2. Also suppose that $\mcA' \in \mbS_{2i}$ and $H' \subseteq \Aut(\mcA')$ are as described above.
Then $p = p(H) = 2i$ and $p' = p(H') = 2i$. 
By parts~(ii) and~(iii) of 
Proposition~\ref{comparisson of structures with different support-substructures},
we have~(\ref{the proportion not having exactly i orbits of cardinality 2 approaches 0}) if
$p(H) - q(H) < p(H') - q(H')$. By assumption we have $p(H) - q(H) = 2i - i = i$.
By assumption, $H'$ has no fixed points, so $H'$ has at most $i$ orbits.
As we also assume that $H'$ does not have $i$ orbits, it follows that
$H'$ has $i'$ orbits for some $i' < i$ and we get
$p(H') - q(H') = 2i - i' > i = p(H) - q(H)$, 
so~(\ref{the proportion not having exactly i orbits of cardinality 2 approaches 0}) follows
from Proposition~\ref{comparisson of structures with different support-substructures}.

We must now prove that there are $\mcA \in \mbS_{2i}$ without fixed point and
a subgroup $H \subseteq \Aut(\mcA)$ without fixed point such that $H$ has 
exactly $i$ orbits. But this holds if we let the interpretation of every relation symbol be
empty (so $\Aut(\mcA) = Sym_{2i}$) and let $H$ the permutation group on $[2i]$ 
with only one nontrivial permutation and this one takes $\alpha$ to $2\alpha$ for every 
$\alpha \in [i]$.

(ii) Suppose that $\mcA \in \mbS_{2i+1}$ has no fixed point and that $H$ is a subgroup
of $\Aut(\mcA)$ without fixed points. Then $p(H) = 2i+1$.
For the same reasons as in part~(i) we only need to show that (subject to the constraint $p(H) = 2i+1$)
$p(H) - q(H)$ is minimal if and only if $H$ has exactly $i$ orbits.
As $H$ has no fixed point it has at most $i$ orbits.
Hence $p(H) - q(H) \geq 2i+1 - i = i+1$ and $p(H) - q(H) = i+1$ if and only if 
$H$ has exactly $i$ orbits. It now suffices to prove that there are $\mcA \in \mbS_{2i+1}$
without fixed point and a subgroup $H \subseteq \Aut(\mcA)$ without fixed point
such that $H$ has exactly $i$ orbits.
If $i = 1$ and we let the interpretation of every relation symbol be empty, then this clearly holds.
So suppose that $i > 1$. Let $B = [2i-2]$ and  $C = \{2i-1, 2i, 2i+1\}$.
Let the interpretation of every relation symbol be empty and let $H \subseteq \Aut(\mcA)$
be the group $H_1 \times H_2$, where
$H_1$ has only one trivial permutation and this one sends $\alpha$ to $2\alpha$ for every $\alpha \in [i-1]$
and fixes every $\alpha \in C$,
every $\alpha \in B$ is a fixed point of $H_2$ and $H_2 \uhrc C$ is the symmetric group of $C$.
Then $\Aut(\mcA) \cong \mbbZ_2 \times Sym_3$ and $\mcA$ has exactly $i$ orbits.
\hfill $\square$

\section{Proof of Theorem~\ref{characterisation of automorphism groups when r=2}}
\label{proof of first main theorem}

\noindent
In this section we prove Theorem~\ref{characterisation of automorphism groups when r=2}.
Lemmas~\ref{typical groups for 2i when r=2}
and~\ref{typical groups for 2i+1 when r=2} 
may be of some interest in themselves.
Throughout this section we assume that
\begin{itemize}{\em
\item[] $r$ is the maximal arity and $r=2$,
\item[] $k$ is the number of $r$-ary relation symbols and
\item[] $l$ is the number of $(r-1)$-ary relation symbols,
}\end{itemize}
although the assumption that $r=2$ is restated in the
results.

\begin{lem}\label{typical groups for 2i when r=2}
Suppose that $i \geq 1$ and $r = 2$.
For almost every $\mcM \in \mbS_n(\spt^* = 2i)$, $\Aut(\mcM) \cong (\mbbZ_2)^t$
for some $t \in \{1, \ldots, i\}$.
Moreover, for every $t \in \{1, \ldots, i\}$ there is a rational number $0 < a_t \leq 1$ such that
\[ \lim_{n\to\infty} 
\frac{\big| \{\mcM \in \mbS_n(\spt^* = 2i) : \Aut(\mcM) \cong (\mbbZ_2)^t\} \big|}
{\big| \mbS_n(\spt^* = 2i) \big|} \ = \ a_t,
\]
and if $i > 1$ then $a_t < 1$.
\end{lem}

\noindent
{\bf Proof.}
By Lemma~\ref{typical number of nontrivial orbits of structures with bounded support},
for almost every $\mcM \in \mbS_n(\spt^* = 2i)$, $\Aut(\mcM) \uhrc \Spt^*(\mcM)$ has 
$i$ orbits, each one of cardinality 2.
For every $\mcM \in  \mbS_n(\spt^* = 2i)$ such that $\Aut(\mcM) \uhrc \Spt^*(\mcM)$ has 
$i$ orbits and every $f \in \Aut(\mcM)$, $f^2$ is the identity.
Hence, for almost every $\mcM \in \mbS_n(\spt^* = 2i)$ there is $t \in \{1, \ldots, i\}$
such that $\Aut(\mcM) \cong (\mbbZ_2)^t$.

By Lemma~\ref{spt-star equals a union of S(A,H)},
there are $\mcA_1, \ldots, \mcA_m \in \mbS_{2i}$ without fixed points and for each
$i = 1, \ldots, m$ a number $l_i$ and subgroups $H_{i,1}, \ldots, H_{i,l_i}$ of $\Aut(\mcA_i)$
without fixed points such that 
\[\mbS_n(\spt^* = 2i) \ = \ \bigcup_{i=1}^m \bigcup_{j=1}^{l_i} \mbS_n(\mcA_i,H_{i,j})\] 
for each sufficiently large $n$.
Recall Lemma~\ref{automorhism group almost surely has same orbits as H}.
Fix $1 \leq t \leq i$. 
Let $\mcA'_i$, $i = 1, \ldots, m$ and $H'_{i,j}$, $j = 1, \ldots, l'_i$, enumerate
all pairs $(\mcA_i, H_{i,j})$ such that $H_{i,j} \cong (\mbbZ_2)^t$
and the proportion of $\mcM \in \mbS_n(\mcA_i, H_{i,j})$ such that $\Aut(\mcM) \cong (\mbbZ_2)^t$
converges to 1.
Now it suffices to prove that 
\[
\frac{\Big| \bigcup_{i=1}^{m'} \bigcup_{j=1}^{l'_i} \mbS_n(\mcA'_i, H'_{i,j}) \Big|}
{\Big| \bigcup_{i=1}^{m} \bigcup_{j=1}^{l_i} \mbS_n(\mcA_i, H_{i,j}) \Big|}
\]
converges to a rational number  as $n \to \infty$. 
But this follows from Proposition~\ref{general quotient proposition}.
Part~(ii)(b) of Proposition~\ref{comparisson of structures with different support-substructures}
guarantees that the limit is larger than 0 if $i > 1$.
\hfill $\square$

\begin{lem}\label{typical groups for 2i+1 when r=2}
Suppose that $i \geq 1$ and $r = 2$.
(i) For almost every $\mcM \in \mbS_n(\spt^* = 3)$, $\Aut(\mcM) \cong \mbbZ_3$
or $\Aut(\mcM) \cong  Sym_3$. 
Moreover, for each one of these groups, call it $G$, there is a rational number $0 < a_G < 1$ such that
\[ \lim_{n\to\infty} 
\frac{\big| \{\mcM \in \mbS_n(\spt^* = 2i+1) : \Aut(\mcM) \cong G\} \big|}
{\big| \mbS_n(\spt^* = 2i+1) \big|} \ = \ a_G.
\]
(ii) Suppose that $i > 1$.
For almost every $\mcM \in \mbS_n(\spt^* = 2i+1)$, $\Aut(\mcM) \cong (\mbbZ_2)^t \times \mbbZ_3$
or $\Aut(\mcM) \cong (\mbbZ_2)^t \times Sym_3$
for some $t \in \{1, \ldots, i-1\}$.
Moreover, for each one of these groups, call it $G$, there is a rational number $0 < a_G < 1$ such that
\[ \lim_{n\to\infty} 
\frac{\big| \{\mcM \in \mbS_n(\spt^* = 2i+1) : \Aut(\mcM) \cong G\} \big|}
{\big| \mbS_n(\spt^* = 2i+1) \big|} \ = \ a_G.
\]
\end{lem}

\noindent
{\bf Proof.}
The first claim of part~(i) is immediate because a permutation group without fixed points
on a set of cardinality 3 must be isomorphic to either $\mbbZ_3$ (if no nonidentity permutation
has a fixed point) or $Sym_3$.
The second claim of part~(i) is proved in the same way as the
second claim of Lemma~\ref{typical groups for 2i when r=2}, with the help of
Propositions~\ref{general quotient proposition}
and~\ref{comparisson of structures with different support-substructures}
and
Lemma~\ref{automorhism group almost surely has same orbits as H}.

Now we prove part~(ii), so suppose that $i > 1$.
By Lemma~\ref{typical number of nontrivial orbits of structures with bounded support},
for almost every $\mcM \in \mbS(\spt^* = 2i+1)$, $\Aut(\mcM) \uhrc \Spt^*(\mcM)$ has 
$i-1$ orbits, say $O_1, \ldots, O_{i-1}$, of cardinality 2 and one orbit $O_i$ of cardinality 3.
Hence, for the first statement of~(ii),
it suffices to prove that for each $\mcM \in \mbS_n(\spt^* = 2i+1)$ with $i-1$ orbits
$O_1, \ldots, O_{i-1}$, of cardinality 2 and one orbit $O_i$ of cardinality 3,
$\Aut(\mcM) \cong (\mbbZ_2)^t \times \mbbZ_3$
or $\Aut(\mcM) \cong (\mbbZ_2)^t \times Sym_3$
for some $t \in \{1, \ldots, i-1\}$.
The second statement of part~(ii) is proved in the same way as the
second statement of part~(i) (and the
second statement of Lemma~\ref{typical groups for 2i when r=2}).

With the given assumptions we have
\[\Aut(\mcM) \uhrc (O_1 \cup \ldots \cup O_{i-1}) \ \cong \ (\mbbZ_2)^t\]
for some $t \geq 1$, because for every $f \in \Aut(\mcM) \uhrc (O_1 \cup \ldots \cup O_{i-1})$,
$f^2$ is the identity.
The next step is to show that 
$\Aut(\mcM) \uhrc \Spt^*(\mcM)$ is the direct product of 
$\Aut(\mcM) \uhrc (O_1 \cup \ldots \cup O_{i-1})$ and
$\Aut(\mcM) \uhrc O_i$,
since it follows from the case $i = 1$ that $\Aut(\mcM) \uhrc O_i \cong \mbbZ_3$ or 
$\Aut \uhrc O_i \cong Sym_3$.

Let $O_i = \{a,b,c\}$.
There is $f \in \Aut(\mcM) \uhrc \Spt^*(\mcM)$ such that $f(a) = b$ and 
$g \in \Aut(\mcM) \uhrc \Spt^*(\mcM)$ such that $g(b) = c$.
If $f(c) = c$ and $g(a) = a$ then $fg$ has no fixed point in $O_i$. 
Otherwise either $f$ or $g$ has no fixed point in $O_i$.
So under all circumstances there exists $f \in \Aut(\mcM) \uhrc \Spt^*(\mcM)$ which has no fixed point in $O_i$.
Since $|O_j| = 2$ for every $j \in \{1, \ldots, i-1\}$ it follows that every
$d \in O_1 \cup \ldots \cup O_{i-1}$ is a fixed point of $f^2$.
Take any $j \in \{1, \ldots, i-1\}$ and let $O_j = \{d,e\}$, so both $d$ and $e$ are fixed points of $f^2$.
Since there is $h \in \Aut(\mcM) \uhrc \Spt^*(\mcM)$ such that $h(d) = e$ (and $h(e) = d$)
it follows, using $f$ and $h$, 
that $O_j \times O_i$ is an orbit of $\Aut(\mcM) \uhrc \Spt^*(\mcM)$ on $\Spt^*(\mcM) \times \Spt^*(\mcM)$.
This holds for every $j \in \{1, \ldots, i-1\}$, and therefore
\[\Aut(\mcM) \uhrc \Spt^*(\mcM) \ \cong \ 
\Aut(\mcM) \uhrc (O_1 \cup \ldots \cup O_{i-1}) \times \Aut(\mcM) \uhrc O_i.\]
Hence, for either $\mcG = \mbbZ_3$ or $\mcG = Sym_3$, and some $t \in \{1, \ldots, i-1\}$,
$\Aut(\mcM) \uhrc \Spt^*(\mcM) \cong (\mbbZ_2)^t \times \mcG$, and clearly the same holds
with $\Aut(\mcM)$ in place of $\Aut(\mcM) \uhrc \Spt^*(\mcM)$.
\hfill $\square$
\\

\noindent
The next corollary is an immediate consequence of 
Lemmas~\ref{typical groups for 2i when r=2}
and~\ref{typical groups for 2i+1 when r=2}.

\begin{cor}\label{almost every structure with support m has an automorphism with support m, case r=2}
Let $m \geq 2$.
Almost every $\mcM \in \mbS(\spt^* = m)$ has an automorphism whose support has cardinality $m$.
\end{cor}

\begin{lem}\label{comparison of 2i+1 over 2i for r=2}
Suppose that $i \geq 1$ and $r = 2$. Then
\[
\lim_{n\to\infty}
\frac{\big| \mbS_n(\spt^* = 2i+1) \big|}{\big| \mbS_n(\spt^* = 2i) \big|}
\ = \ 0.
\]
\end{lem}

\noindent
{\bf Proof.}
By Lemma~\ref{typical number of nontrivial orbits of structures with bounded support},
for almost every $\mcM \in \mbS_n(\spt^*  = 2i)$,
$H = \Aut(\mcM) \uhrc \Spt^*(\mcM)$ has exactly $i$ orbits,
and for almost every $\mcM' \in \mbS_n(\spt^*  = 2i+1)$,
$H' = \Aut(\mcM') \uhrc \Spt^*(\mcM')$ has exactly $i$ orbits.
For such $H$ and $H'$ we have 
\[p(H) - q(H) \ = \ i \ < \ i+1 = p(H') - q(H'),\]
so if $\mcA = \mcM \uhrc \Spt^*(\mcM)$ and $\mcA' = \mcM' \uhrc \Spt^*(\mcM')$
(and $\mcM$ and $\mcM'$ are as above), then
Proposition~\ref{comparisson of structures with different support-substructures}
implies that
$\big|S_n(\mcA',H')\big| \big/ \big|S_n(\mcA,H)\big| \to 0$ as $n \to \infty$.
The lemma follows from this because,
by Lemma~\ref{spt-star equals a union of S(A,H)}, 
each one of $\mbS(\spt^* = 2i)$ and $\mbS(\spt^* = 2i+1)$
is a union of finitely many sets of the form $\mbS(\mcA,H)$.
\hfill $\square$

\begin{lem}\label{comparison of 2i+2 over 2i+1 for r=2}
Suppose that $i \geq 1$ and $r = 2$. Then
\[
\lim_{n\to\infty}
\frac{\big| \mbS_n(\spt^* = 2i+1) \big|}{\big| \mbS_n(\spt^* = 2i+2) \big|}
\ = \ 0.
\]
\end{lem}

\noindent
{\bf Proof.}
By Lemma~\ref{typical number of nontrivial orbits of structures with bounded support},
for almost every $\mcM \in \mbS_n(\spt^*  = 2i+2)$,
$H = \Aut(\mcM) \uhrc \Spt^*(\mcM)$ has exactly $i+1$ orbits,
and for almost every $\mcM' \in \mbS_n(\spt^*  = 2i+1)$,
$H' = \Aut(\mcM') \uhrc \Spt^*(\mcM')$ has exactly $i$ orbits.
It follows that
\[p(H) - q(H) \ = 2i + 2 - (i+1) \ = \ i+1 \ = \ 2i+1 - i \ = \ p(H') - q(H')\]
and
\[p(H) \ = \ 2i + 2 \ > \  2i + 1 \ = \ p(H').\]
So if $\mcM$ and $\mcM'$ are as above, $\mcA = \mcM \uhrc \Spt^*(\mcM)$ and
$\mcA' = \mcM' \uhrc \Spt^*(\mcM')$, then
Proposition~\ref{comparisson of structures with different support-substructures} implies
that $\big| S_n(\mcA',H') \big| \big/ \big| S_n(\mcA,H) \big| \to 0$ as $n \to \infty$.
The lemma follows because each one of $\mbS(\spt^* = 2i+1)$ and $\mbS(\spt^* = 2i+2)$
is a union of finitely many sets of the form $\mbS(\mcA,H)$.
\hfill $\square$

\begin{lem}\label{comparison of m+2 over m for r=2}
Suppose that $r = 2$ and either $m = 0$ or $m \geq 2$. Then
\[
\lim_{n\to\infty}
\frac{\big| \mbS_n(\spt^* = m+2) \big|}{\big| \mbS_n(\spt^* = m) \big|}
\ = \ 0.
\]
\end{lem}

\noindent
{\bf Proof.}
The case $m = 0$ follows from the fact that almost all $\mcM \in \mbS$ are rigid \cite{Fag77}.
Now suppose that $m \geq 2$.
By Lemma~\ref{typical number of nontrivial orbits of structures with bounded support},
for almost every $\mcM \in \mbS_n(\spt^*  = m)$,
$H = \Aut(\mcM) \uhrc \Spt^*(\mcM)$ has exactly $\lfloor \frac{m}{2} \rfloor$ orbits,
and for almost every $\mcM' \in \mbS_n(\spt^*  = m+2)$,
$H' = \Aut(\mcM') \uhrc \Spt^*(\mcM')$ has exactly 
$\lfloor \frac{m+2}{2} \rfloor = \lfloor \frac{m}{2} \rfloor + 1$ orbits.
Since
\[ p(H) - q(H) \ = m - \Big\lfloor \frac{m}{2} \Big\rfloor \ < \ 
m - \Big\lfloor \frac{m}{2} \Big\rfloor + 1 \ = \ 
p(H') - q(H'),\]
it follows that if $\mcM$ and $\mcM'$ are as above, $\mcA = \mcM \uhrc \Spt^*(\mcM)$
and $\mcA' = \mcM' \uhrc \Spt^*(\mcM')$, then 
Proposition~\ref{comparisson of structures with different support-substructures}
implies that
$\big| S_n(\mcA', H') \big| \big/ \big| S_n(\mcA,H) \big| \to 0$ as $n \to \infty$,
which in turn implies the lemma (just as in the proofs of the preceeding two lemmas).
\hfill $\square$

\begin{lem}\label{comparison between m and T when r=2}
Suppose that $r = 2$.
Also assume that $m = 0$ or $m \geq 2$ and that $T > m$ and $T \geq 2$.
Let $m' = m$ if $m$ is even and $m' = m+1$ otherwise.
Then
\[\lim_{n\to\infty} \frac{\big| \mbS_n(\spt^* = m') \big|}{\big| \mbS_n(m \leq \spt^* \leq T) \big|}
\ = \ 1.\]
\end{lem}

\noindent
{\bf Proof.}
The case when $m = 0$ follows from \cite{Fag77}, 
so suppose that $m \geq 2$.
If $T = m+1$ then the result follows from 
Lemmas~\ref{comparison of 2i+1 over 2i for r=2}
and~\ref{comparison of 2i+2 over 2i+1 for r=2}.
Now suppose that $m \geq 2$ and $T \geq m+2$.
For each $i \in \{m+2, \ldots, T\}$ we have, by Lemma~\ref{comparison of m+2 over m for r=2},
\[
\frac{\big| \mbS_n(\spt^* = i) \big|}{\big| \mbS_n(m \leq \spt^* \leq T) \big|} \ \leq \ 
\frac{\big| \mbS_n(\spt^* = i) \big|}{\big| \mbS_n(\spt^* = i-2)\big|} \ \to \ 0
\]
as $n \to \infty$. From this it follows that
\[
\lim_{n\to\infty} 
\frac{\big| \mbS_n(\spt^* = m) \cup \mbS_n(\spt^* = m+1) \big|}{\big| \mbS_n(m \leq \spt^* \leq T) \big|}
\ = \ 1.
\]
The lemma now follows from
Lemmas~\ref{comparison of 2i+1 over 2i for r=2}
and~\ref{comparison of 2i+2 over 2i+1 for r=2}.
\hfill $\square$

\begin{lem}\label{almost all structures with spt at least m have spt* = m, case r=2}
Suppose that $r = 2$ and $m \geq 2$.
Let $m' = m$ if $m$ is even and $m' = m+1$ otherwise.
Then
\[\lim_{n\to\infty} \frac{\big| \mbS_n(\spt \geq m) \cap \mbS_n(\spt^* = m') \big|}{\big| \mbS_n(\spt \geq m) \big|}
\ = \ 
\lim_{n\to\infty} \frac{\big| \mbS_n(\spt^* = m') \big|}{\big| \mbS_n(\spt^* \geq m) \big|}
\ = \ 1.\]
\end{lem}

\noindent
{\bf Proof.}
Let $m \geq 2$.
Proposition~\ref{bound on support}
says that there is $T > m$ such that
\[
\lim_{n\to\infty} 
\frac{\big| \mbS_n(m \leq \spt^* \leq T) \big|}{\big| \mbS_n(\spt^* \geq m) \big|} \ = \ 
\lim_{n\to\infty} 
\frac{\big| \mbS_n(\spt \geq m) \ \cap \ \mbS_n(\spt^* \leq T) \big|}{\big| \mbS_n(\spt \geq m) \big|} \ = \ 1.
\]
By Corollary~\ref{almost every structure with support m has an automorphism with support m, case r=2}
it suffices to prove that
\[
\lim_{n\to\infty} \frac{\big| \mbS_n(\spt^* = m') \big|}{\big| \mbS_n(m \leq \spt^* \leq T) \big|} \ = \ 1,
\]
but this follows from Lemma~\ref{comparison between m and T when r=2}.
\hfill $\square$
\\

\noindent
We get Theorem~\ref{characterisation of automorphism groups when r=2}
by combining Lemmas~\ref{typical groups for 2i when r=2}
and~\ref{almost all structures with spt at least m have spt* = m, case r=2}.

\section{Proof of Theorem~\ref{characterisation of automorphism groups when r>2}}
\label{proof of second main theorem}

\noindent
Theorem~\ref{characterisation of automorphism groups when r>2}
is proved in this section, but
Lemmas~\ref{the number of nontrivial orbits on pairs}
and~\ref{the almost sure automorphism group when spt*=2i+1} 
may be of interest in themselves.
{\em The symbols $r$, $k$ and $l$ have the same meaning in this section as in the previous};
see the beginning of it.
{\em But in this section we assume that $r \geq 3$.}

\begin{lem}\label{the number of nontrivial orbits on pairs}
Suppose that $r \geq 3$ and $i \geq 1$. 
For almost all $\mcM \in \mbS(\spt^* = 2i)$,
$\Aut(\mcM) \uhrc \Spt^*(\mcM)$ has exactly
$2i^2$ orbits on $\Spt^*(\mcM) \times \Spt^*(\mcM)$, so every orbit on
$\Spt^*(\mcM) \times \Spt^*(\mcM)$ has cardinality 2.
\end{lem}

\noindent
{\bf Proof.}
First note that if $H$ is the permutation group on 
$\Omega = \{v_1, \ldots, v_i, w_1, \ldots, w_i\}$
whose only nontrivial
permutation sends $v_j$ to $w_j$ for every $j$,
then $H$ has $i$ orbits on $\Omega$ and $2i^2$ orbits on $\Omega \times \Omega$,
because every orbit on $\Omega \times \Omega$ has cardinality 2.
Hence $s(H) = 2i^2$.
Moreover, for every permutation group on $\Omega$ without fixed points,
the number of orbits on $\Omega \times \Omega$ cannot exceed $(2i)^2/2 = 2i^2$.
So if $H$ is as described then $s(H)$ is maximal among permutation groups on a set of cardinality $2i$.
We also have $p(H) - q(H) = i$ which is minimal among permutation groups without any fixed point
on a set of cardinality $2i$.
Let $\mcA$ be any structure without fixed point with universe $A = \Omega$ such that $H$ 
is a subgroup of $\Aut(\mcA)$.
For example, let the interpretation of every relation symbol be empty. 
Suppose that $\mcA'$ is a structure with universe of cardinality $2i$ and without any fixed point 
and suppose, moreover, that $H'$ is a subgroup of $\Aut(\mcA')$ such that $H'$ has no
fixed point and either $q(H') < i$ or $s(H') < 2i^2$.
By Proposition~\ref{comparisson of structures with different support-substructures},
$\big| \mbS_n(\mcA',H') \big| \Big/ \big| \mbS_n(\mcA,H) \big| \ \to \ 0$ as $n \to \infty$.
By Lemma~\ref{automorhism group almost surely has same orbits as H},
almost all $\mcM \in \mbS(\mcA, H)$ have the property that 
the number of orbits of $\Aut(\mcM) \uhrc \Spt^*(\mcM)$ on $\Spt^*(\mcM)$ is $q(H) = i$
and the number of orbits of $\Aut(\mcM) \uhrc \Spt^*(\mcM)$ on $\Spt^*(\mcM) \times \Spt^*(\mcM)$ is $s(H) = 2i^2$.
Now the lemma follows, because $\mbS(\spt^* = 2i)$ is a union of finitely many sets of the form
$\mbS(\mcA,H)$ where the universe of $\mcA$ has cardinality $2i$, $\mcA$ has no fixed point
and $H$ is subgroup of $\Aut(\mcA)$ without fixed point.
\hfill $\square$

\begin{lem}\label{the almost sure automorphism group when spt*=2i}
Suppose that $r \geq 3$ and $i \geq 1$. 
For almost all $\mcM \in \mbS(\spt^* = 2i)$, $\Aut(\mcM) \cong \mbbZ_2$.
\end{lem}

\noindent
{\bf Proof.}
Since $\Aut(\mcM) \uhrc \Spt(\mcM) \cong \Aut(\mcM)$ it suffices to prove that
for almost all $\mcM \in \mbS(\spt^* = 2i)$, $\Aut(\mcM) \uhrc \Spt(\mcM) \cong \mbbZ_2$.
By Lemmas~\ref{typical number of nontrivial orbits of structures with bounded support}
and~\ref{the number of nontrivial orbits on pairs}
it suffices to prove that if $H$ is a permutation group on $[2i]$ such that
every orbit of $H$ on $[2i]$ has cardinality 2 and every orbit of
$H$ on $[2i] \times [2i]$ has cardinality 2, then $H \cong \mbbZ_2$.
This is obvious if $i = 1$, so for the rest of the proof we assume that $i \geq 2$.

So suppose that $H$ is a permutation group on $[2i]$ such that
every orbit of $H$ on $[2i]$ has cardinality 2 and 
every orbit of $H$ on $[2i] \times [2i]$ has cardinality 2.
We first prove an auxilliary claim.
\\

\noindent
{\em Claim.} If $a$ and $b$ belong to different orbits of $H$ on $[2i]$ and
$f \in H$ is not the identity, then $\{f(a), f(b)\} \cap \{a,b\} = \es$.
\\

\noindent
{\em Proof.}
Suppose for a contradiction that the claim does not hold.
Then there are orbits $\{a,c\},\{b,d\} \subseteq [2i]$ and   
a permutation $f \in H$ such that $f(a) = c$ and $f(b) = b$.
Then $f(d) = d$ and as $\{b,d\}$ is an orbit there is
$g \in H$ such that $g(b) = d$ and $g(d) = b$.
If $g(a) = a$ then $\{a,c\} \times \{b,d\}$ is an orbit
of $H$ on $[2i] \times [2i]$, contradicting the assumption that all
orbits on $[2i] \times [2i]$ have cardinality 2.
Hence $g(a) = c$ and $g(c) = a$.
Then $gf(a) = a$ and $gf(b) = d$ and again, by the use of $f$, $gf$ and
compositions of them, it follows that
$\{a,c\} \times \{b,d\}$ is an orbit, contradicting our assumption.
\hfill $\square$
\\

\noindent
Now we prove that if $f \in H$ is not the identity, then $f$ has no fixed point.
Suppose, for a contradiction, that $f \in H$ is not the identity and has a fixed point $a$.
As the orbit to which $a$ belongs, say $\{a,c\}$, has cardinality 2 and we assume that $i \geq 2$
it follows that there is $b \in [2i] \setminus \{a,c\}$ such that $f(b) \neq b$. 
Then we have $a = f(a) \in \{f(a), f(b)\} \cap \{a,b\}$, contradicting the claim.

Next, we prove that $H$ has a unique nonidentity permutation from which it
follows that $H \cong \mbbZ_2$.
So suppose for a contradiciton that $f, g \in H$ are 
nonidentity permutations and $f(a) \neq g(a)$ for some $a$.
Then $a, f(a)$ and $g(a)$ belong to the same orbit.
Since neither $f$ nor $g$ has any fixed point, as we proved above, some
orbit of $H$ on $[2i]$ contains at least three elements, contradicting our assumption.
\hfill $\square$
\\

\noindent
The next result deals only with permutation groups and is independent of
the ingredients from formal logic such as relation symbols and their interpretations.

\begin{lem}\label{auxilliary lemma for Z-2 times Z-3}
Suppose that $i \geq 2$.
Let $H$ be a permutation group without fixed points on $[2i+1]$ such
that $H$ has exactly $i-1$ orbits of cardinality 2, exactly one orbit of cardinality
3 and no other orbits. 
If $s(H)$ is maximal among all $H$ subject to the given constraints,
then $H \cong \mbbZ_2 \times \mbbZ_3$ and $s(H) = 2i^2 - 2i + 3$.
\end{lem}

\noindent
{\bf Proof.}
Suppose that $H$ is a permutation group without fixed points on $[2i+1]$ such
that $H$ has exactly $i-1$ orbits of cardinality 2, exactly one orbit of cardinality
3 and no other orbits.
Let $O_1, \ldots, O_{i-1}$ be the orbits with cardinality 2 and let $O_i$ be the
orbit with cardinalty 3. 
Let $\Omega = O_1 \cup \ldots \cup O_{i-1}$

We first show that if $H \uhrc \Omega \cong \mbbZ_2$,
$H \uhrc O_i \cong \mbbZ_3$ and $H \cong (H \uhrc \Omega) \times (H \uhrc O_i)$, then
$s(H) = 2i^2 - 2i + 3$.
So suppose that $H \uhrc \Omega \cong \mbbZ_2$.
Then $H \uhrc \Omega$ has exactly $i-1$ orbits on $\Omega$, each one of cardinality 2,
and $H \uhrc \Omega$ has exactly $2(i-1)^2$ orbits on $\Omega \times \Omega$.
Now suppose that $H \uhrc O_i \cong \mbbZ_3$.
Then it is easy to see that no $f \in H \uhrc O_i$ other than the identity 
has a fixed point in $O_i$ and therefore
$H \uhrc O_i$ has exactly 3 orbits on $O_i \times O_i$.
Suppose, in addition to previous assumptions and conclusions, that
$H \cong (H \uhrc \Omega) \times (H \uhrc O_i)$.
Then it easily follows that for every $j = 1, \ldots, i-1$, $O_j \times O_i$ and $O_i \times O_j$ are orbits
of $H$ on $[2i+1] \times [2i+1]$.
Hence, the number of orbits of $H$ on $[2i+1] \times [2i+1]$ which contain $(a,b)$ such
that $a \in \Omega$ and $b \in O_i$, or vice versa, is $2(i-1)$.
Altogether, we get
\[s(H) \ = \ 2(i-1)^2 + 3 + 2(i-1) \ = \ 2i^2 - 2i + 3.\]
We now show that if $s(H)$ is maximal among all $H$ subject to the given constraints in the lemma,
then $H \cong \mbbZ_2 \times \mbbZ_3$. This will conclude the proof.

So suppose that $s(H)$ is maximal among all permutation groups on $[2i+1]$ without fixed points
and with exactly $i-1$ orbits of cardinality 2, exactly one orbit of cardinality
3 and no other orbits. 
As before, let $O_1, \ldots, O_{i-1}$ be the orbits with cardinality 2, let $O_i$ be the
orbit with cardinalty 3 and let $\Omega = O_1 \cup \ldots \cup O_{i-1}$.
By the same argument as in the proof of Lemma~\ref{typical groups for 2i+1 when r=2}
there exists $f \in H$ without any fixed point in $O_i$, and then, just as in that proof,
it follows that for every $j = 1, \ldots, i-1$, $O_j \times O_i$ and $O_i \times O_j$
are orbits of $H$ on $[2i+1] \times [2i+1]$.
Hence the number of orbits of $H$ on $[2i+1] \times [2i+1]$ that contain a pair $(a,b)$ such that
$a \in \Omega$ and $b \in O_i$, or vice versa, is at most $2(i-1)$.
The number of orbits of $H$ on $[2i+1] \times [2i+1]$ that contain a pair $(a,b)$ where
$a,b \in \Omega$ is at most $(2(i-1))^2/2 = 2(i-1)^2$, because every orbit has at least two members.
It is easy to see that number of orbits of $H$ on $[2i+1] \times [2i+1]$ that contain a pair $(a,b)$ where 
$a,b \in O_i$ is at most 3 (one orbit containing $(a,a)$ where $a \in O_i$, one orbit containing
$(a,b)$ for some distinct $a,b \in O_i$ and one orbit containing $(b,a)$).
This means that
\[s(H) \ \leq \ 2(i-1) + 2(i-1)^2 + 3 \ = \ 2i^2 - 2i + 3.\]
By the assumption that $s(H)$ is maximal and since the value
$2i^2 - 2i + 3$ can be reached, as shown above, we get
$s(H) = 2i^2 - 2i + 3$.
From the argument above it follows that $s(H)$ cannot be maximal unless
$H \uhrc \Omega$ has a maximal number of orbits on $\Omega \times \Omega$.
Hence $H \uhrc \Omega$ must have the maximal possible number of orbits on $\Omega \times \Omega$
which is $(2(i-1))^2/2 = 2(i-1)^2$  and consequently
every orbit of $H \uhrc \Omega$ on $\Omega \times \Omega$ has cardinality 2.
By the argument in the proof of Lemma~\ref{the almost sure automorphism group when spt*=2i}
it follows that $H \uhrc \Omega \cong \mbbZ_2$.

We have seen that $H \uhrc O_i$ can have at most 3 orbits on $O_i \times O_i$.
Also it is easy to see that $H \uhrc O_i$ has 3 orbits on $O_i \times O_i$
if and only if for any distinct $a,b \in O_i$, $(a,b)$ and $(b,a)$ belong to different orbits.
Moreover, if for any distinct $a,b \in O_i$, $(a,b)$ and $(b,a)$ belong to
different orbits, then no $f \in H \uhrc O_i$ has order 2,
so $H \uhrc O_i \cong \mbbZ_3$.

By the same argument as in the proof of Lemma~\ref{typical groups for 2i+1 when r=2},
using only the assumptions about the orbits of $H$ on $\Omega$, it follows 
that $H \cong (H \uhrc \Omega) \times (H \uhr O_i) \cong \mbbZ_2 \times \mbbZ_3$.
\hfill $\square$
\\

\noindent
Recall that by Definition~\ref{definition of beta}:
\[
\beta(x,y,z) \ = \ k\binom{r}{2}x^2 \ - \ kr(r-1)xy \ - \ l(r-1)x \ + \ l(r-1)y \ + \ k\binom{r}{2}z.
\]
Also remember that $k$ is the number of $r$-ary relation symbols and $l$ is the number
of $(r-1)$-ary relation symbols.

\begin{lem}\label{the almost sure automorphism group when spt*=2i+1}
Suppose that $r \geq 3$.\\
(i) For almost all $\mcM \in \mbS_n(\spt^* = 3)$, $\Aut(\mcM) \cong \mbbZ_3$.\\
(ii) If $i \geq 2$ then for almost all $\mcM \in \mbS(\spt^* = 2i+1)$, 
$\Aut(\mcM) \cong \mbbZ_2 \times \mbbZ_3$
and $s\big(\Aut(\mcM) \uhrc \Spt^*(\mcM)\big) = 2i^2 - 2i + 3$.
\end{lem}

\noindent
{\bf Proof.}
We start with part~(ii), so suppose that $i \geq 2$. 
Suppose that $\mcA \in \mbS_{2i+1}$ has no fixed point and suppose that $H$ is a subgroup
of $\Aut(\mcA)$ without fixed point. Note that $p(H) = 2i+1$. We have seen, 
in the proof of 
Lemma~\ref{typical number of nontrivial orbits of structures with bounded support}~(ii), 
that $p(H) - q(H)$ is minimal when $q(H) = i$ 
(under the assumption that $H$ acts on a set of cardinality $2i+1$ and has no fixed points), 
which implies that $H$ has $i-1$ orbits of cardinality 2 and one orbit
of cardinality 3.
Also, recall the definition of
$\beta(x,y,z)$ in Corollary~\ref{estimate of the number of structures for given p,q,s when r>2}.
Observe that if $p = p(H) = 2i+1$, $q = q(H) = i$ and $s = s(H)$, then
\[\beta(p,q,s) \ = \ k\binom{r}{2}(2i+1)^2 \ - \ kr(r-1)(2i+1)i \ - \ l(r-1)(2i+1) 
\ + \ l(r-1)i \ + \ k\binom{r}{2}s,\]
where $r, k, l$ and $i$ are fixed parameters.
So under the assumptions that $p(H) = 2i+1$ and $q(H) = i$, $\beta(p,q,s)$ is maximised when
$s = s(H)$ is maximised. 
From Proposition~\ref{comparisson of structures with different support-substructures}~(iii)
and the fact that $\mbS(\spt^* = 2i+1)$ is a union of finitely
many sets of the form $\mbS(\mcA, H)$, where $\mcA \in \mbS_{2i+1}$, 
$\mcA$ has no fixed point and $H$ is a subgroup of $\Aut(\mcA)$ without any fixed point, 
it follows that almost every $\mcM \in \mbS(\spt^* = 2i+1)$ has the following
properties:
$H = \Aut(\mcM) \uhrc \Spt^*(\mcM)$ has exactly $i$ orbits ($i-1$ of cardinality 2 and one of cardinality 3)
and $s(H)$ is maximal among all permutation groups on $[2i+1]$
with $i$ orbits and without a fixed point.
From Lemma~\ref{auxilliary lemma for Z-2 times Z-3} it now follows that for
almost every $\mcM \in \mbS(\spt^* = 2i+1)$, 
$\Aut(\mcM) \cong \mbbZ_2 \times \mbbZ_3$ and 
$s\big(\Aut(\mcM) \uhrc \Spt^*(\mcM)\big) = 2i^2 - 2i + 3$.

Now we consider part~(i). 
If $\mcM \in \mbS(\spt^* = 3)$ and $H = \Aut(\mcM) \uhrc \Spt^*(\mcM)$ then $p(H) = 3$
and $q(H) = 1$. 
So the question is what $H = \Aut(\mcM) \uhrc \Spt^*(\mcM)$ looks like when $s(H)$,
the number of orbits on $\Spt^*(\mcM) \times \Spt^*(\mcM)$, is maximised.
It is easy to see that
$s(H) \leq 3$, and $s(H) = 3$ if and only if $H \cong \mbbZ_3$.
(We argued similarly in Lemma~\ref{auxilliary lemma for Z-2 times Z-3}.)
\hfill $\square$
\\

\noindent
From Lemmas~\ref{the almost sure automorphism group when spt*=2i}
and~\ref{the almost sure automorphism group when spt*=2i+1} we get the following:

\begin{cor}\label{almost every structure with support m has an automorphism with support m, case r>2}
Suppose that $r \geq 3$ and $m \geq 2$.
Almost every $\mcM \in \mbS(\spt^* = m)$ has an automorphism whose support has cardinality $m$.
\end{cor}

\begin{lem}\label{comparison of 2i+1 over 2i for r>2}
Suppose that $i \geq 1$ and $r \geq 3$. Then
\[
\lim_{n\to\infty}
\frac{\big| \mbS_n(\spt^* = 2i+1) \big|}{\big| \mbS_n(\spt^* = 2i) \big|}
\ = \ 0.
\]
\end{lem}

\noindent
{\bf Proof.}
Exactly as the proof of Lemma~\ref{comparison of 2i+1 over 2i for r=2}.
\hfill $\square$

\begin{lem}\label{auxilliary lemma for beta with different inputs}
\[\beta(2i+2, \ i+1, \ 2(i+1)^2) \ - \ \beta(2i+1, \ i, \ 2i^2-2i+3) \ = \ 
2k\binom{r}{2}(2i - 1).\]
\end{lem}

\noindent
{\bf Proof.}
Straightforward, but long, calculation.
\hfill $\square$

\begin{lem}\label{comparison of 2i+2 over 2i+1 for r>2}
If $r \geq 3$ and $i \geq 1$ then
\[
\lim_{n\to\infty} \frac{\big| \mbS_n(\spt^* = 2i+1) \big|}{\big| \mbS_n(\spt^* = 2i+2) \big|}
\ = \ 0. 
\]
\end{lem}

\noindent
{\bf Proof.}
By Lemmas~\ref{typical number of nontrivial orbits of structures with bounded support}
and~\ref{the number of nontrivial orbits on pairs}, 
for almost all $\mcM \in \mbS_n(\spt^* = 2i+2)$, if $H = \Aut(\mcM) \uhrc \Spt^*(\mcM)$
then $p = p(H) = 2i+2$, $q = q(H) = i+1$ and $s = s(H) = 2(i+1)^2$.
By Lemmas~\ref{typical number of nontrivial orbits of structures with bounded support} 
and~\ref{the almost sure automorphism group when spt*=2i+1}, 
for almost all $\mcM' \in \mbS_n(\spt^* = 2i+1)$, if $H' = \Aut(\mcM') \uhrc \Spt^*(\mcM)$
then $p' = p(H') = 2i+1$, $q' = q(H') = i$ and $s' = s(H') = 2i^2 - 2i + 3$.
For such $H$ and $H'$ we have 
\[p - q  \ = \  i+1 \ = \ p' - q' \]
and by Lemma~\ref{auxilliary lemma for beta with different inputs} we also have
\[\beta(p,q,s) \ > \ \beta(p',q',s'),\]
so Lemma~\ref{comparison of 2i+2 over 2i+1 for r>2} follows from
Proposition~\ref{comparisson of structures with different support-substructures}~(iii).
\hfill $\square$

\begin{lem}\label{comparison of m+2 over m for r>2}
Suppose that $r \geq 3$ and either $m = 0$ or $m \geq 2$. Then
\[
\lim_{n\to\infty}
\frac{\big| \mbS_n(\spt^* = m+2) \big|}{\big| \mbS_n(\spt^* = m) \big|}
\ = \ 0.
\]
\end{lem}

\noindent
{\bf Proof.}
Exactly as the proof of Lemma~\ref{comparison of m+2 over m for r=2}.
\hfill $\square$

\begin{lem}\label{comparison between m and T when r>2}
Suppose that $r \geq 3$ and suppose that $m = 0$ or $m \geq 2$.
Let $m' = m$ if $m$ is even and $m' = m+1$ otherwise. Then
For every integer $T$ such that $T > m$ and $T \geq 2$,
\[
\lim_{n\to\infty} \frac{\big| \mbS_n(\spt^* = m') \big|}{\big| \mbS_n(m \leq \spt^* \leq T) \big|}
\ = \ 1.
\]
\end{lem}

\noindent
{\bf Proof.}
As the proof of Lemma~\ref{comparison between m and T when r=2},
but now using Lemmas~\ref{comparison of 2i+1 over 2i for r>2},
~\ref{comparison of 2i+2 over 2i+1 for r>2}
and~\ref{comparison of m+2 over m for r>2}.
\hfill $\square$

\begin{lem}\label{almost all structures with spt at least m have spt* = m, case r>2}
Suppose that $r \geq 3$ and $m \geq 2$.
Let $m' = m$ if $m$ is even and $m' = m+1$ otherwise.
Then
\[\lim_{n\to\infty} \frac{\big| \mbS_n(\spt \geq m) \cap \mbS_n(\spt^* = m') \big|}{\big| \mbS_n(\spt \geq m) \big|}
\ = \ 
\lim_{n\to\infty} \frac{\big| \mbS_n(\spt^* = m') \big|}{\big| \mbS_n(\spt^* \geq m) \big|}
\ = \ 1.\]
\end{lem}

\noindent
{\bf Proof.}
Like the proof of Lemma~\ref{almost all structures with spt at least m have spt* = m, case r=2},
but now using 
Corollary~\ref{almost every structure with support m has an automorphism with support m, case r>2}
and
Lemma~\ref{comparison between m and T when r>2}.
\hfill $\square$
\\

\noindent
By combining Lemmas~\ref{the almost sure automorphism group when spt*=2i}
and~\ref{almost all structures with spt at least m have spt* = m, case r>2}
we get Theorem~\ref{characterisation of automorphism groups when r>2}.

\end{document}